\documentclass[a4paper]{amsart}

\usepackage{amsmath,amstext,amssymb,mathrsfs,amscd,amsthm,indentfirst, bbm}
\usepackage{amsfonts, mathtools, afterpage}
\usepackage[dvipsnames,svgnames,x11names,hyperref]{xcolor}
\definecolor{darkred}{RGB}{139,0,0}
\definecolor{darkblue}{RGB}{0,0,139}
\definecolor{darkgreen}{RGB}{0,100,0}
\definecolor{darkyellow}{RGB}{100,100,0}
\usepackage{ifthen, nicefrac, scalerel}
\usepackage[shortlabels]{enumitem}

\usepackage{tikz, tikz-cd}
\usetikzlibrary{matrix,calc,positioning,arrows,decorations.pathreplacing,decorations.markings,patterns}
\tikzset{->-/.style={decoration={
			markings,
			mark=at position #1 with {\arrow{>}}},postaction={decorate}}}
\tikzset{-<-/.style={decoration={
					markings,
					mark=at position #1 with {\arrow{<}}},postaction={decorate}}}

\usepackage{booktabs}
\usepackage{verbatim}

\usepackage{graphicx}

\usepackage{amsthm}

\theoremstyle{plain}

\newtheorem{theorem}{Theorem}[section]

\theoremstyle{definition}

\newcommand{\bQ}{\mathbb{Q}}
\newcommand{\bR}{\mathbb{R}}
\newcommand{\bS}{\mathbb{S}}

\newcommand{\bZ}{\mathbb{Z}}

\newcommand\lra{\longrightarrow}

\newcommand\Diff{\mathrm{Diff}}
\newcommand\Tor{\mathrm{Tor}}

\newcommand\Sp{\mathrm{Sp}}
\newcommand\Homeo{\mathrm{Homeo}}
\newcommand\hAut{\mathrm{hAut}}
\newcommand\Top{\mathrm{Top}}

\newcommand\OO{\mathrm{O}}
\newcommand\SO{\mathrm{SO}}
\newcommand\CC{\mathrm{C}}
\newcommand\fr{\mathrm{fr}}

\newcommand\Emb{\mathrm{Emb}}

\newcommand\Bun{\mathrm{Bun}}

\newcommand\hocolim{\operatorname*{hocolim}}
\newcommand\holim{\operatorname*{holim}}

\newcommand{\map}{\mathrm{map}}

\mathchardef\ordinarycolon\mathcode`\:
\mathcode`\:=\string"8000
\begingroup \catcode`\:=\active
  \gdef:{\mathrel{\mathop\ordinarycolon}}
\endgroup

\numberwithin{equation}{section}

\usepackage[hypertexnames=false]{hyperref}

\title{Diffeomorphisms of discs}

\author{Oscar Randal-Williams}
\thanks{O.\ Randal-Williams was supported by the ERC under the European Union’s Horizon 2020 research and innovation programme (grant agreement No.\ 756444), and by a Philip Leverhulme Prize from the Leverhulme Trust.}
\email{o.randal-williams@dpmms.cam.ac.uk}
\address{Centre for Mathematical Sciences\\
Wilberforce Road\\
Cambridge CB3 0WB\\
UK}

\begin{document}
\begin{abstract}
I describe what is currently known, for $d \geq 5$, about the rational homotopy type of the group of diffeomorphisms of the $d$-disc relative to its boundary, and the closely related group of homeomorphisms of $d$-dimensional Euclidean space.
\end{abstract}
\maketitle

\section{Introduction}

I will be concerned with the homotopy types of the topological groups of diffeomorphisms $\Diff(\bR^d)$ and homeomorphisms $\Homeo(\bR^d)$ of Euclidean space, and of diffeomorphisms $\Diff_\partial(D^d)$ and homeomorphisms $\Homeo_\partial(D^d)$ of closed discs fixing the boundary pointwise (equivalently: compactly-supported diffeomorphisms and homeomorphisms of Euclidean space). By scaling outwards the group $\Diff(\bR^d)$ deformation retracts to the subgroup of linear diffeomorphisms, and thence to the subgroup $\OO(d)$ of orthogonal diffeomorphisms of $\bR^d$. By scaling inwards (the ``Alexander trick'') the group $\Homeo_\partial(D^d)$ contracts to a point.

In contrast the groups $\Homeo(\bR^d)$ and $\Diff_\partial(D^d)$ have much more mysterious homotopy types. As long as $d \neq 4$ all four groups are related by smoothing theory \cite[Essay V]{kirbysiebenmann}, which provides a homotopy equivalence
$$\tfrac{\Homeo_\partial(D^d)}{\Diff_\partial(D^d)} \simeq \Omega^d_0 \left( \tfrac{\Homeo(\bR^d)}{\Diff(\bR^d)}\right)$$
and so, incorporating the above and writing $\Top(d) := \Homeo(\bR^d)$, provides a homotopy equivalence (the ``Morlet equivalence'')
$$B\Diff_\partial(D^d) \simeq \Omega^d_0 \tfrac{\Top(d)}{\OO(d)}.$$
Taking the homotopy type of $\OO(d)$ as given, the Morlet equivalence shows that understanding the homotopy types of $\Top(d)$ and $B\Diff_\partial(D^d)$ are more or less equivalent. The latter is independently interesting as it classifies smooth $D^d$-bundles $\pi : E \to B$ trivialised near the boundary, and this perspective offers a useful way to study it: it is the perspective I usually adopt. For any manifold $M$ of dimension $d\neq 4$ smoothing theory identifies $\tfrac{\Homeo_\partial(M)}{\Diff_\partial(M)}$ with certain path components of a space of sections $\Gamma_\partial(\mathrm{Fr}(TM) \times_{\OO(d)} \tfrac{\Top(d)}{\OO(d)} \to M)$, so these homotopy types furthermore describe the difference between diffeomorphisms and homeomorphisms of all $d$-manifolds. It is therefore an important goal of geometric topology to investigate these homotopy types.

In this essay I will describe what is known about the \emph{rational} homotopy type of $B\Diff_\partial(D^d)$, and some recent techniques which are allowing us to say more about it. Along the way the influence of Michael Weiss will be seen at every turn, and it is a pleasure to acknowledge and celebrate his many profound contributions to this subject.

\section{Some phenomena}

I will first describe the classical approach to calculating $\pi_*(B\Diff_\partial(D^d))\otimes\bQ$, which describes it completely in the so-called pseudoisotopy stable range, and then explain two more recent results which indicate the existence of new phenomena outside of this range: the work of Watanabe on configuration space integrals, and the work of Weiss on unstable topological Pontrjagin classes.

\subsection{Pseudoisotopy and algebraic $K$-theory}\label{sec:concordances}

The topological group of smooth pseudoisotopies 
$$\CC(M) := \{f : M \times [0,1] \overset{\text{diffeo}}\lra M \times [0,1] \, | \, f \text{ fixes } \partial M \times [0,1] \cup M \times \{0\} \text{ pointwise}\}$$
of a manifold $M$ participates in a fibre sequence
\begin{equation}\label{eq:ConcFibSeq}
\Diff_\partial(M \times [0,1]) \lra \CC(M) \xrightarrow{f \mapsto f\vert_{M \times \{1\}}} \Diff_\partial(M),
\end{equation}
and so measures to what extent diffeomorphisms of $M \times [0,1]$ may be represented as loops of diffeomorphisms of $M$. There is a stabilisation map
$$\sigma_M : \CC(M) \lra \CC(M \times [0,1])$$
morally induced by crossing with the interval, but technically slightly more involved. The (smooth) \emph{pseudoisotopy stable range} $\Phi(d)$ is the minimum of the connectivities of $\sigma_M$ taken over all manifolds $M$ of dimension $\geq d$,  and it is a deep theorem of Igusa \cite{igusastab} that $\Phi(d) \geq \min(\tfrac{d-7}{2}, \tfrac{d-4}{3})$. The stabilisation $\mathscr{C}(M) := \hocolim_n \CC(M \times [0,1]^n)$ may be promoted to a homotopy-invariant functor from the category of spaces to the category of infinite loop spaces. The stable parameterised $h$-cobordism theorem \cite{waldhausenstable} 
relates the functor $\mathscr{C}(-)$ to Waldhausen's algebraic $K$-theory of spaces---or to the $K$-theory of ring spectra---by the fibration sequence (with cosection) of infinite loop spaces
$$B\mathscr{C}(M) \lra Q(M_+) \lra \Omega^\infty \mathrm{K}(\mathbb{S}[\Omega M]).$$

In particular, for $M=D^d$ the rational homotopy equivalence $\mathrm{K}(\mathbb{S}) \to \mathrm{K}(\bZ)$ and Borel's calculation \cite{borelstable} 
\begin{equation*}
K_i(\bZ)\otimes\bQ = \begin{cases}
\bQ & i = 0, 5, 9, 13, 17, 21, \ldots\\
0 & \text{else}
\end{cases}
\end{equation*}
determines $\pi_*(\CC(D^d))\otimes\bQ$ in the pseudoisotopy stable range as being a copy of $\bQ$ in each degree $\equiv 3 \mod 4$.

A further piece of structure on $\CC(M)$ is the \emph{pseudoisotopy involution}. Writing $\tau$ for the reflection of $[0,1]$ at $\tfrac{1}{2}$, this is given by
$$f \longmapsto \bar{f} =  (f\vert_{M \times \{1\}} \times [0,1])^{-1} \circ (M \times \tau) \circ f \circ (M \times \tau),$$
and there are compatible involutions on the fibration sequence \eqref{eq:ConcFibSeq} given by inversion on $\Diff_\partial(M)$, and by conjugating by the reflection $M \times \tau$ on $\Diff_\partial(M \times [0,1])$. By analysing this involution, as well as the other involution on $\CC(M \times [0,1])$ induced by $\CC(M \times \tau)$ and their compatibility with $\sigma_M$, it can be shown (see \cite[Section 6.5]{igusatorsion} for a nice discussion) that
\begin{equation}\label{eq:FH}
\begin{aligned}
\pi_i(B\Diff_\partial(D^{2n}))\otimes\bQ &= 0\\
\pi_i(B\Diff_\partial(D^{2n+1}))\otimes\bQ &= \begin{cases}
\bQ & i \equiv 0 \mod 4\\
0 & \text{else}
\end{cases}
\end{aligned}
\end{equation}
in the pseudoisotopy stable range. This calculation was first obtained by Farrell and Hsiang \cite{farrellhsiang}, though by somewhat different means.

\subsection{Configuration space integrals}\label{sec:CSI}

Kontsevich \cite{KontsevichECM} proposed a method to produce invariants of smooth \emph{vertically framed} $D^{d}$-bundles $\pi : E \to B$ trivialised near the boundary, by forming (certain compactifications of) the fibrewise configuration spaces and integrating suitably chosen differential forms along these associated configuration space bundles. The combinatorics of these forms are organised in terms of graphs, and the result is a chain map $\mathrm{GC}_{d}^2[-d]^\vee \otimes_\bQ \bR \to \Omega^\bullet(B)$ from a certain graph complex to the de Rham complex of the base, so the homology of this graph complex yields invariants of the original bundle. Up to regrading, the chain complexes $\mathrm{GC}_{d}^2$ only depend on the parity of $d$, and they split $\mathrm{GC}_{d}^2 = \bigoplus_{g} \mathrm{GC}_{d}^{2, g\text{-loop}}$ as a sum of subcomplexes of fixed loop order. See Willwacher's contribution to the 2018 ICM for an introduction to these objects.

\vspace{3ex}

\noindent\textbf{The work of Watanabe.}
The detailed investigation of this construction has been taken up by Watanabe, firstly in \cite{watanabe1, watanabe2, WatanabeIIerr} for odd-dimensional discs as Kontsevich proposed, and latterly \cite{WatanabeS4, WatanabePC} for even-dimensional discs too. Write $B\Diff^\fr_\partial(D^{d})$ for the space classifying smooth vertically framed $D^{d}$-bundles trivialised near the boundary. (The analogue of the Morlet equivalence in this setting has the form $B\Diff^\fr_\partial(D^{d}) \simeq \Omega^{d} \Top(d)$---up to a small correction of path-components which I shall ignore---so studying smooth framed disc bundles has an even closer connection to homeomorphisms of Euclidean space.) Kontsevich's construction in particular gives characteristic classes
$$\xi_r \in H^{r \cdot (d-3)}(B\Diff^\fr_\partial(D^{d}) ; \mathcal{A}_r^{(-1)^d}),$$
where $\mathcal{A}_r^{+}$ and $\mathcal{A}_r^{-}$ are real vector spaces spanned by connected trivalent graphs of loop order $r+1$ (equipped with certain orientation data which I shall neglect), modulo the IHX relation (and modulo certain signs when changing orientation data: these signs depend on the supercript $+$ and $-$). This vector space arises as the lowest nontrivial homology of $\mathrm{GC}_{d}^{2, (r+1)\text{-loop}} \otimes_\bQ \bR$: the differential is given by summing over splitting vertices, so all trivalent graphs are automatically cycles, and the IHX relation arises from the three ways to split a 4-valent vertex. For small values of $r$ the dimension of $\mathcal{A}_r^{-}$ has been calculated to be $1,1,1,2,2,3,4,5,6,8,9$ for $r=1,2,\ldots, 11$ and the dimension of $\mathcal{A}_r^{+}$ has been calculated to be $0,1,0,0,1,0,0,0,1$ for $r=1,2,\ldots, 9$

Watanabe's results in this direction (\cite[Theorem 3.1]{watanabe2} taking into account the improvement in \cite{WatanabeIIerr}, and \cite{WatanabeS4} taking into account the improvement in \cite{WatanabePC}) is that as long as $d \geq 4$ the evaluation map
$$\xi_r : \pi_{r \cdot (d-3)}(B\Diff^\fr_\partial(D^{d}))\otimes\bR \lra \mathcal{A}_r^{(-1)^d}$$
is surjective. In fact his result is somewhat more precise: he constructs for each trivalent graph $\Gamma$ a more-or-less explicit framed $D^{d}$-bundle over a sphere which is sent by the map $\xi_r$ to the class of $\Gamma$.

This does not directly tell us about $B\Diff_\partial(D^{d})$ because of the framing data, but the difference is easily understood. Forgetting framings defines a homotopy fibre sequence
$$\Omega^{d} \OO(d) \lra B\Diff^\fr_\partial(D^{d}) \lra B\Diff_\partial(D^{d}),$$
and it is not hard to calculate
$$\pi_*(\Omega^{d} \OO(d))\otimes \bQ = \bigoplus_{\substack{k \geq 3,\\ k \equiv 2d+1 \mod 4}} \bQ[d-k].
$$
Thus as long as $d$ is even, or $d$ is odd and $r > 1$, one still has the lower bound 
$$\dim_\bQ \pi_{r \cdot (d-3)}(B\Diff_\partial(D^{d}))\otimes\bQ \geq \dim_\bR \mathcal{A}_r^{(-1)^d}.$$ 
It is worth pausing at this point to emphasise how remarkable it is that Watanabe's results apply for $d=4$, and imply for example that $\pi_2(B\Diff_\partial(D^{4}))\otimes\bQ \neq 0$.

However when $d$ is odd---say $d=2n+1$---and $r=1$, the composition
$$\bR \cong \pi_{2n-2}(\Omega^{2n+1} \OO(2n+1))\otimes\bR \lra \pi_{2n-2}(B\Diff^\fr_\partial(D^{2n+1}))\otimes\bR \xrightarrow{\xi_1} \mathcal{A}_1^\mathrm{odd} \cong \bR$$
might be nontrivial, and in fact it is. Watanabe addresses this difficulty in his earlier work \cite{watanabe1}, by constructing an integral refinement of $\xi_1$, and playing off this integrality against the non-integrality of the topological Pontrjagin classes: his conclusion is that as long as certain arithmetic conditions (\cite[Corollary 2]{watanabe1}, \cite[Corollary 3.5]{watanabe2}) involving Bernoulli numbers and the orders of stable homotopy groups of spheres is satisfied, then one may still conclude that $\pi_{2n-2}(B\Diff_\partial(D^{2n+1}))\otimes\bQ \neq 0$. He verified this computationally for all odd $n \leq 399$. 

\vspace{3ex}

\noindent\textbf{Automorphisms of the little discs operads.}
Configuration space integrals are constructed from (suitable compactifications of) all the ordered configuration spaces $\mathrm{Conf}_k(D^d)$, and all the natural maps between them, applied fibrewise to a (vertically framed) $D^d$-bundle. There is another way to encode (the homotopy types of) these configuration spaces and the natural maps between them, namely as the little $d$-discs operad $E_d$. There is a topological version of the framed little $d$-discs operad, which means that in a suitable homotopical sense the group $\Top(d)$ acts on $E_d$, giving a map
$$B\Top(d) \lra B\hAut(E_d),$$
where the latter is the classifying space of the $E_1$-algebra of derived automorphisms of the little $d$-discs operad. Looping this$(d+1)$ times gives a map
$$B\Diff^\fr_\partial(D^{d}) \simeq \Omega^{d} \Top(d) \lra \Omega^d \hAut(E_d).$$
(This corresponds \cite{boavidaweisscats} to applying the embedding calculus of Goodwillie and Weiss \cite{weissembeddings, goodwillieweiss} to framed self-embeddings of $D^d$ relative to the boundary, though that point of view is not necessary for this discussion.)

The derived automorphisms of the \emph{rationalised} little $d$-discs operad $E_d^\bQ$ have been studied by Fresse, Turchin, and Willwacher \cite{FTW}, who for $d \geq 3$ give an identification $\pi_i(\hAut(E_d^\bQ)) = H_i(\mathrm{GC}_{d}^2)$. Combined with the above this gives a map
$$\pi_i(B\Diff_\partial^\fr(D^d))\otimes\bQ \lra H_{i+d}(\mathrm{GC}_{d}^2)$$
and it is difficult to imagine that this is given by anything other than evaluation of Kontsevich's invariant, but as far as I know the connection between this point of view and configuration space integrals has not yet been made precise. Assuming for now that this is so, Watanabe's results show that this map hits those graph homology classes represented by trivalent graphs.

\subsection{Pontrjagin--Weiss classes}\label{sec:PWclasses}

It follows from the work of Sullivan and of Kirby and Siebenmann that the homotopy fibre $\Top/\OO$ of the map $B\OO \to B\Top$ has finite homotopy groups, and therefore that
$$H^*(B\Top;\bQ) = \bQ[p_1, p_2, p_3, \ldots],$$
a polynomial ring on certain classes $p_i$ of degree $4i$ which pull back to the Pontrjagin classes on $B\OO$: these are the topological Pontrjagin classes. By pulling back along the stabilisation map $B\Top(d) \to B\Top$ they are defined for all $\bR^d$-bundles.

For real vector bundles of dimension $2n$, and so universally in the cohomology of $B\OO(2n)$, the definition of Pontrjagin classes in terms of Chern classes of the complexification immediately gives that
\begin{align}
p_i &= 0 \text{ for } i > n \label{eq:WeissQuestion1}\\
p_n &= e^2 \label{eq:WeissQuestion2}
\end{align}
where $e$ denotes the Euler class. The Euler class only depends on the underlying spherical fibration of the vector bundle, obtained by removing the zero-section. An $\bR^d$-bundle also has an associated spherical fibration (by removing any section), and hence also has an Euler class: as both the Euler and Pontrjagin classes are defined on $B\Top(2n)$, one may then ask about the validity of the identities \eqref{eq:WeissQuestion1} and \eqref{eq:WeissQuestion2} there.

\vspace{3ex}

\noindent\textbf{The work of Weiss.}
Reis and Weiss \cite{ReisWeiss2016} had proposed an elaborate strategy for establishing these identities, but in a spectacular turnaround Weiss \cite{WeissDalian} then showed that these identities are in fact generally \emph{false} for $\bR^d$-bundles. I will comment further on his strategy in Section \ref{sec:WeissSeq}, as its philosophy is fundamental to all the results in Section \ref{sec:HtyType}.

To say more precisely what Weiss proved, consider the fibration sequence
$$\tfrac{\Top(2n)}{\OO(2n)} \lra B\OO(2n) \lra B\Top(2n).$$
The rational cohomology classes $p_n - e^2$ and $p_i$ for $i > n$ are defined on $B\Top(2n)$ and are canonically trivial on $B\OO(2n)$, and hence yield (pre-)transgressed cohomology classes $(p_n-e^2)^\tau$ and $p_i^\tau$ on $\tfrac{\Top(2n)}{\OO(2n)}$. Weiss showed \cite[Section 6]{WeissDalian} that for many $n$ and $i \geq n$ (he shows that $n \geq 83$ and $i < \tfrac{9n}{4}-11$, or $n \geq 59$ and $i < \tfrac{7n}{4}$ will do) these evaluate nontrivially against $\pi_{4i-1}(\tfrac{\Top(2n)}{\OO(2n)})$: this certainly implies that the corresponding $p_n-e^2$ and $p_i$ are nontrivial in the cohomology of $B\Top(2n)$, but is stronger. Translated to diffeomorphisms groups of discs via the Morlet equivalence, Weiss' result shows that the map
$$\pi_{4i-2n-1}(B\Diff_\partial(D^{2n})) \cong \pi_{4i-2n-1}\left(\Omega^{2n}_0\tfrac{\Top(2n)}{\OO(2n)}\right) \cong \pi_{4i-1}\left(\tfrac{\Top(2n)}{\OO(2n)}\right) \xrightarrow{(p_n-e^2)^\tau \text{ or } p_i^\tau} \bQ,$$
is nontrivial for many $n$ and $i \geq n$. I will call an element of $\pi_{4i-2n-1}(B\Diff_\partial(D^{2n}))$ a \emph{Pontrjagin--Weiss class} if it is detected by such maps.

\vspace{3ex}

\noindent\textbf{Odd dimensions.}
On $B\OO(2n+1)$ the Pontrjagin classes still satisfy $p_i=0$ for $i> n$, which is the analogue of \eqref{eq:WeissQuestion1}, so there are (pre-)transgressed classes $p_i^\tau$ on $\tfrac{\Top(2n+1)}{\OO(2n+1)}$ for $i > n$ and one may ask about their nontriviality on homotopy groups. As these classes pull back to the classes of the same name on $\tfrac{\Top(2n)}{\OO(2n)}$, this nontriviality follows from Weiss' theorem in many cases.

It seems to be less well known that there is also an analogue of \eqref{eq:WeissQuestion2} in odd dimensions. Namely, there is a characteristic class $E$ of $S^{2n}$-fibrations such that
\begin{align}
p_n &= E \label{eq:WeissQuestion3}
\end{align}
in $H^{4n}(B\OO(2n+1);\bQ)$. It may be defined as follows. Given a $S^{2n}$-fibration $S^{2n} \to X \overset{\pi}\to Y$ with orientation local system $\mathscr{O}$, the self-intersection of the fibrewise diagonal map $\Delta : X \to X \times_Y X$ defines a \emph{fibrewise Euler class} $e^\mathrm{fw}(\pi) \in H^{2n}(X;\mathscr{O})$, and then $E(\pi) := \tfrac{1}{2} \int_\pi e^\mathrm{fw}(\pi)^3 \in H^{4n}(Y;\bQ)$.)

As $E$ depends only on the underlying spherical fibration it is also defined in $H^*(B\Top(2n+1);\bQ)$, so one can also ask whether the identity \eqref{eq:WeissQuestion3} fails to hold here, or even better whether the cohomology class $(p_n-E)^\tau$ on $\tfrac{\Top(2n+1)}{\OO(2n+1)}$ evaluates nontrivially on $\pi_{4n-1}(\tfrac{\Top(2n+1)}{\OO(2n+1)})$. It can be checked that under $B\Top(2n) \to B\Top(2n+1)$ the class $E$ pulls back to $e^2$, so $(p_n-E)^\tau$ pulls back to $(p_n - e^2)^\tau$ on $\tfrac{\Top(2n)}{\OO(2n)}$, and hence the nontriviality of $(p_n-E)^\tau$ on homotopy groups in many degrees also follows from Weiss' theorem.

\vspace{3ex}

\noindent\textbf{Propagating.}
Formalising the method used above, the stabilisation maps
\[
\begin{tikzcd}
\Omega^{d}_0 \tfrac{\Top(d-1)}{\OO(d-1)} \dar \arrow[r, equals] & \Omega_0 B\Diff_\partial(D^{d-1}) \dar\\
\Omega^{d}_0 \tfrac{\Top(d)}{\OO(d)} \arrow[r, equals]&  B\Diff_\partial(D^{d})
\end{tikzcd}
\]
(known as ``Gromoll  maps" on the diffeomorphism group side) show that if a Pontrjagin--Weiss class exists on $B\Diff_\partial(D^{d-1})$, and the cohomology class detecting it can be defined on $\tfrac{\Top(d)}{\OO(d)}$, then it survives to $B\Diff_\partial(D^{d})$.

\vspace{3ex}

\noindent\textbf{Relation to configuration space integrals.}
Somewhat surprisingly the map
$$\pi_{2n-2}(B\Diff^\fr_\partial(D^{2n+1})) \cong \pi_{4n}(B\Top(2n+1)) \xrightarrow{E} \bQ$$
can be identified with the simplest Kontsevich invariant $\xi_1$ (which is that associated to the $\Theta$-graph) as studied by Watanabe in \cite{watanabe1}. From this point of view Watanabe's argument in that paper shows that $p_n \neq E$ in $H^{4n}(B\Top(2n+1);\bQ)$, so is closely related to Weiss' theorem (but does not imply it!). This is explained in detail in \cite[Appendix B]{KrannichR-WOdd}.

\section{The rational homotopy type of $B\Diff_\partial(D^d)$}\label{sec:HtyType}

The results of the last section give a complete calculation \eqref{eq:FH} of the rational homotopy groups $\pi_*(B\Diff_\partial(D^d))\otimes\bQ$ valid in the pseudoisotopy stable range, but also indicate the existence of various new phenomena outside of this range. These new phenomena start in degrees $\sim d$, and Krannich \cite{KrannichConc} and I \cite{oscarconcordance} had shown that \eqref{eq:FH} is in fact valid in degrees $\lesssim d$. In this section I present two recent results, obtained in collaboration with Kupers and with Krannich, giving detailed information quite far outside of this range, and I then speculate about what they might be indicating.

\subsection{Even-dimensional discs}

Kupers and I \cite{KR-WDisks, KR-WTorelliLie} have investigated the rational homotopy type of $B\Diff_\partial(D^{2n})$. The following is the main result of \cite{KR-WDisks}, incorporating the improvement from \cite[Section 7.1]{KR-WTorelliLie}.

\begin{theorem}[Kupers--R-W]\label{thm:KupersRW}
Let $2n \geq 6$. Then $\pi_j(B\Diff_\partial(D^{2n}))\otimes\bQ=0$ for $j < 2n-1$, and for $j \geq 2n-1$ we have
$$\pi_j(B\Diff_\partial(D^{2n}))\otimes\bQ = \begin{cases}
\bQ & j \equiv 2n-1 \mod 4, j \not\in \bigcup\limits_{r \geq 2} [2r(n-2)-1, 2r(n-1)+1]\\
0 & j \not\equiv 2n-1 \mod 4, j \not\in \bigcup\limits_{r \geq 2} [2r(n-2)-1, 2r(n-1)+1]\\
? & \text{otherwise}.
\end{cases}$$
\end{theorem}
The copies of $\bQ$ in this theorem are generated by Pontrjagin--Weiss classes in the sense of Section \ref{sec:PWclasses}, and the theorem gives a complete calculation in degrees $\leq 4n-10$, as well as in higher degrees outside of the indicated ``bands''. 

It can be cast in a somewhat stronger form by using the Morlet equivalence $B\Diff_\partial(D^{2n}) \simeq \Omega^{2n}_0 \tfrac{\Top(2n)}{\OO(2n)}$ and considering the fibration sequence
$$ \Omega^{2n+1} \tfrac{\Top}{\Top(2n)} \lra \Omega^{2n} \tfrac{\Top(2n)}{\OO(2n)} \lra \Omega^{2n} \tfrac{\Top}{\OO(2n)}.$$
A slight strengthening of the theorem is then that $\pi_*( \Omega^{2n+1} \tfrac{\Top}{\Top(2n)})\otimes\bQ $ is supported in degrees $\bigcup_{r \geq 2} [2r(n-2)-1, 2r(n-1)+1]$; the rational homotopy groups of $\Omega^{2n} \tfrac{\Top}{\OO(2n)}$ are $\bQ$ in every degree $\equiv 2n-1 \mod 4$, and the right-hand map detects the Pontrjagin--Weiss classes.

The result can also be given a little more structure by using the involution on $B\Diff_\partial(D^{2n}) \simeq \Omega^{2n}_0 \tfrac{\Top(2n)}{\OO(2n)}$ induced by conjugation by a reflection of the disc. The terms in the fibration sequence above have compatible involutions, which on $\pi_*(\Omega^{2n}_0 \tfrac{\Top}{\OO(2n)})\otimes\bQ$ acts as $(-1)$, and on $\pi_*(\Omega^{2n+1} \tfrac{\Top}{\Top(2n)})\otimes\bQ$ acts as $(-1)^r$ in the band of degrees $[2r(n-2)-1, 2r(n-1)+1]$ (when such bands overlap this should be regarded as inconclusive). This implies the existence of Pontrjagin--Weiss classes outside of degrees $\bigcup_{r \geq 2, r \text{ odd}} [2r(n-2)-1, 2r(n-1)+1]$.

Finally, as explained in Section \ref{sec:PWclasses}, Pontrjagin--Weiss classes can be propagated from smaller discs to larger ones. The conclusion of this discussion is depicted in Figure \ref{fig:1}. It seems likely that all possible Pontrjagin--Weiss classes already exist in $\pi_*(B\Diff_\partial(D^6))$.

\afterpage{%
\begin{figure}[hbtp]
\centering
	\begin{tikzpicture}
	\begin{scope}[scale=.65]
	
	\def\HH{25} 
	\def\WW{18} 
	\def\HHhalf{13} 
	\def\WWhalf{10} 
	
	\def\AA{-1}
	\def\BB{1}

	\clip (-1,-1) rectangle ({\WW+0.5},{\HH+0.25});
	\draw (-.5,0)--({\WW+.5},0);
	\draw (0,-1) -- (0,{\HH+1.5});

	\foreach \s in {2,4,6,8,10,12,14,16,18,20}
	{
		\fill[black!20!white,opacity=.5] (0.5,{\s*(3-2-0.5)+.5*\AA}) -- (\WW,{\s*\WW+.5*\AA}) -- (\WW,{\s*(\WW+1)+.5*\BB}) -- (0.5,{\s*(3-1-0.5)+.5*\BB});
	}

	\foreach \s in {3,5,7,9,11,13,15,17,19,21}
	{
		\fill[black!20!white,opacity=.8] (0.5,{\s*(3-2-0.5)+.5*\AA}) -- (\WW,{\s*\WW+.5*\AA}) -- (\WW,{\s*(\WW+1)+.5*\BB}) -- (0.5,{\s*(3-1-0.5)+.5*\BB});
	}

	
	  
	\node [fill=white] at (-.5,-.75) {$\nicefrac{*}{d}$};

	\draw [thick,black,dotted] (1,{0.25}) -- (4.5,4/2);
	\draw [thick,black,dotted] (4.5,4/2) -- (\WW, {(2*\WW+3)/6});

	\def\PWdestab{{1,1,2,2,3,3,4,6,5,7,8,8,9,10,12,11,11,11,11,11,11,11,11,11}}
	\foreach \s in {0,...,23}
	{
	\pgfmathsetmacro\PWlen{\PWdestab[\s]};
		\foreach \i in {1,...,\PWlen}
		{
			\node [black] at (\s-\i+2,\s+\i+1.5) {$\bullet$};
		}
	}

	\def\Aplus{{0,1, , ,1, , , ,1}}
	
	\foreach \s in {2,...,9}
	{
					\pgfmathsetmacro\TheNum{\Aplus[\s-1]};
	        \foreach \x in {1,...,\WW}
        {
          \node [black] at (\x,\s*\x+0.5*\s) {$\scaleobj{0.6}{\mathbf{\TheNum}}$};
        }
	}
	
		\begin{scope}
	\foreach \s[evaluate={\si=int(2*\s)}] in {0,...,\HH}
	{
		\draw [dotted] (-.5,\s)--(.25,\s);
		\draw [dotted] (.75,\s) -- ({\WW+.5},\s);
		\node [fill=white] at (-.25,\s) [left] {\tiny $\si$};
	}

	\foreach \s[evaluate={\si=int(2*(\s+2))}] in {1,...,\WW}
	{
		\draw [dotted] (\s,-0.5)--(\s,{\HH+.5});
		\node [fill=white] at (\s,-.5) {\tiny $\si$};
	}
	\end{scope}

	\end{scope}
	\end{tikzpicture}

\caption{Rational homotopy groups of $B\Diff_\partial(D^{2n})$. The calculation is complete in the unshaded region, and $\bullet$ denotes Pontrjagin--Weiss classes. The lightly shaded bands are those on which the reflection acts as $+1$, and Pontrjagin--Weiss classes are still present in these; the darkly-shaded bands are those where the reflection acts as $-1$. Existing copies of $\bullet$ have been propagated downwards along lines of slope $-1$ as in Section \ref{sec:PWclasses}. The numbers denote Watanabe's lower bounds on these groups. The dotted line indicates the Igusa stable range.}\label{fig:1}
\end{figure}
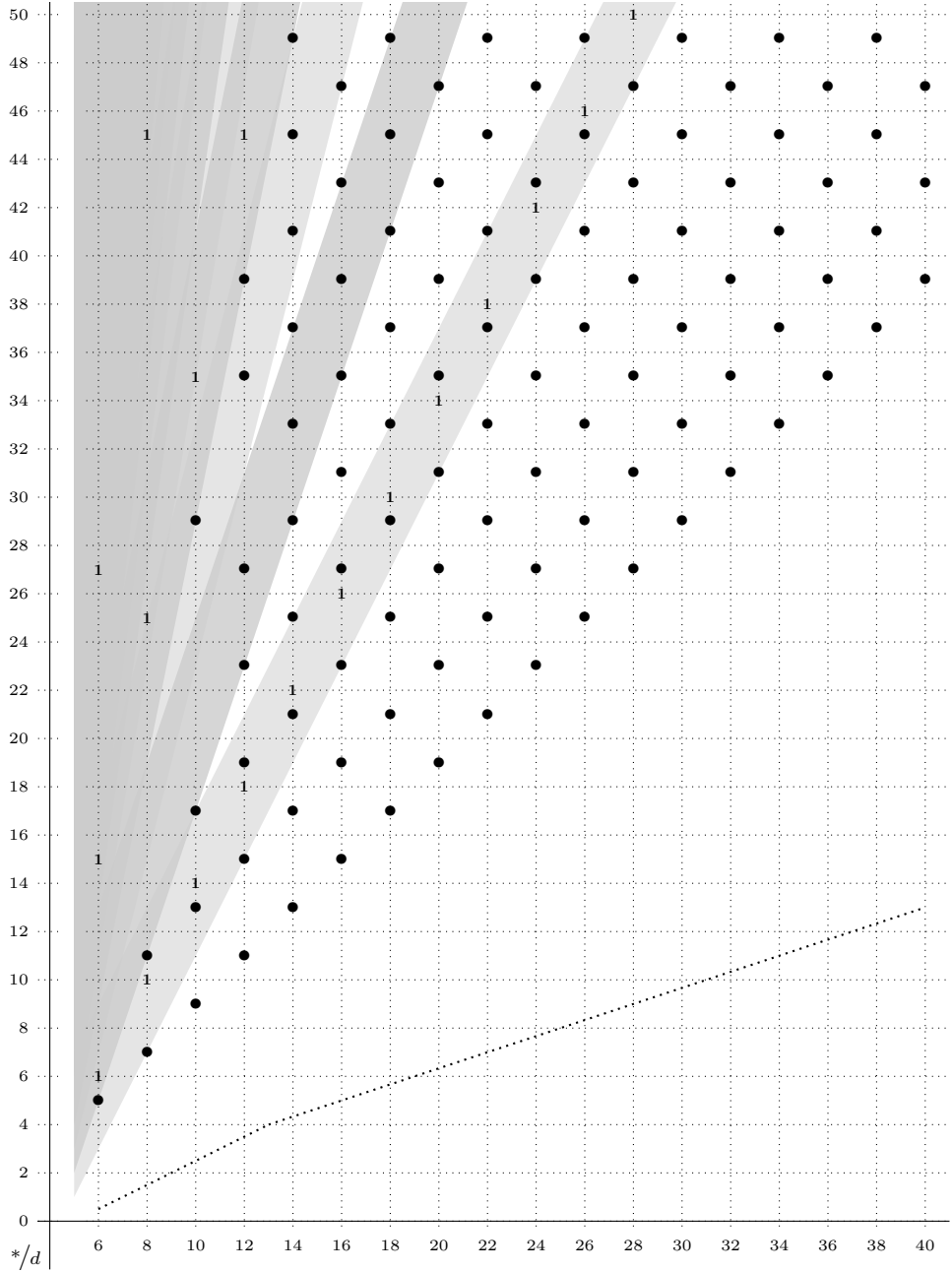
\clearpage
}

\afterpage{%
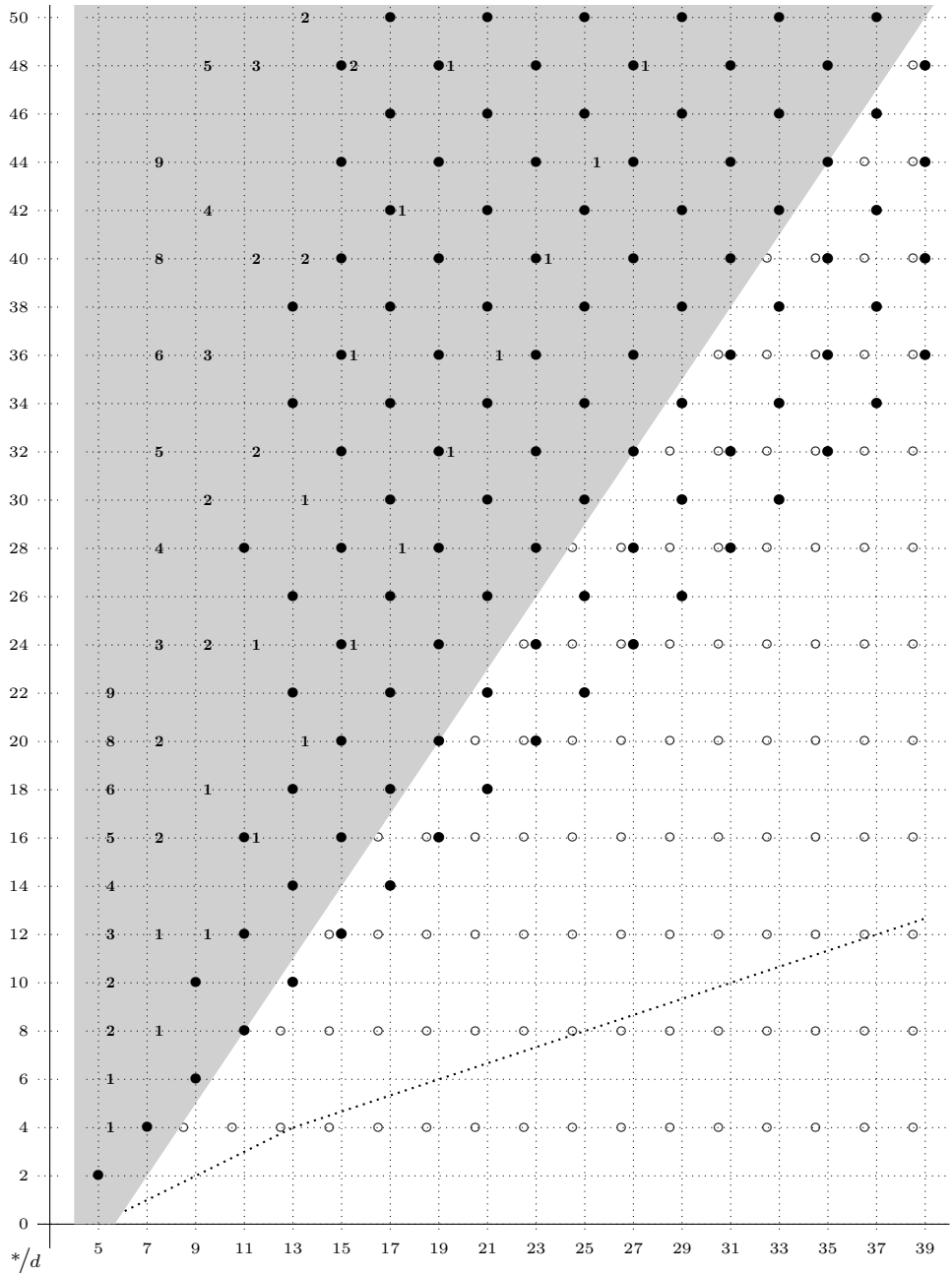
\begin{figure}[hbtp]
\centering

        \begin{tikzpicture}
        \begin{scope}[scale=.65]

        \def\HH{25} 
        \def\WW{18} 
        \def\HHhalf{13} 
        \def\WWhalf{10} 

        \clip (-1,-1) rectangle ({\WW+0.5},{\HH+0.25});
        \draw (-.5,0)--({\WW+.5},0);
        \draw (0,-.5) -- (0,{\HH+1.5});

	\draw [thick,black,dotted] (1,0) -- (5,2);
	\draw [thick,black,dotted] (5,2) -- (\WW, {(2*\WW+2)/6});

        \foreach \x in {1,...,\WW}
        {
                \foreach \k[evaluate={\ki=int(6*\k-1)}] in {1,...,\HHhalf}
                {
								\node [black] at (\x-0.25,{2*\k}) {$\scaleobj{0.9}{\circ}$};
                }
        }

                \foreach \x[evaluate={\xi=int(\x-0.25)}] in {1,...,\WW}
        {
                \foreach \k in {1,...,\HHhalf}
                {
                \node [black] at (\x,{\x-2+2*\k}) {$\bullet$};
                }
        }
				
				\definecolor{mygray}{HTML}{D0D0D0}
				\fill[mygray] (1.333,0) -- (\WW+1,{1.5*\WW-0.5}) -- (\WW+1,\HH+1) -- (0.5,\HH+1) -- (0.5,0);


		\node [black] at (1,{1}) {$\bullet$};
		
		
	\def\PWdestab{{1,1,2,2,3,3,4,6,5,7,8,8,9,10,12,11,11,11,11,11,11,11,11,11}}
	\foreach \s in {0,...,23}
	{
	\pgfmathsetmacro\PWlen{\PWdestab[\s]};
		\foreach \i in {1,...,\PWlen}
		{
			\node [black] at (\s-\i+2+1,\s+\i+1) {$\bullet$};
		}
	}
		
	\def\Aminus{{1,1,1,2,2,3,4,5,6,8,9}}
	
	\foreach \s in {2,...,11}
	{
					\pgfmathsetmacro\TheNum{\Aminus[\s-1]};
	        \foreach \x in {1,...,\WW}
        {
          \node [black] at (\x+0.25,\s*\x) {$\scaleobj{0.6}{\mathbf{\TheNum}}$};
        }
	}

	

	
	

				
				



   \begin{scope}

        \foreach \y[evaluate={\yi=int(2*\y)}] in {0,...,\HH}
        {
										\draw [dotted] (-.5,\y)--(.25,\y);
										\draw [dotted] (.75,\y) -- ({\WW+.5},\y);
                \node [fill=white] at (-.25,\y) [left] {\tiny $\yi$};
        }       

        \foreach \x[evaluate={\xi=int(2*\x+3)}] in {1,...,\WW}
        {
                \draw [dotted] (\x,-0.5)--(\x,{\HH+.5});
                \node [fill=white] at (\x,-.5) {\tiny $\xi$};
        }
     \end{scope}

        \node [fill=white] at (-.5,-.75) {$\nicefrac{*}{d}$};

        \end{scope}
     \end{tikzpicture}

\caption{Rational homotopy groups of $B\Diff_\partial(D^{2n+1})$. The calculation is complete in the unshaded region: $\bullet$ denotes Pontrjagin--Weiss classes, $\circ$ denotes algebraic $K$-theory classes. In the shaded region we have indicated existing Pontrjagin--Weiss classes, and the numbers denote Watanabe's lower bounds on these groups.}\label{fig:2}
\end{figure}
\clearpage
}

\subsection{Odd-dimensional discs} Krannich and I \cite{KrannichR-WOdd} have investigated the rational homotopy type of $B\Diff_\partial(D^{2n+1})$.

\begin{theorem}[Krannich--R-W]\label{thm:KrannichRW}
In degrees $j \leq 3n-8$ we have
\begin{equation*}
\pi_j(B\Diff_\partial(D^{2n+1}))\otimes\bQ = K_{j+1}(\bZ)\otimes\bQ  \oplus \begin{cases}
\bQ & j \equiv 2n-2 \mod 4, j \geq 2n-2\\
0 & \text{otherwise}.
\end{cases}
\end{equation*}
\end{theorem}

The first term is the rational algebraic $K$-theory of the integers, extending the classes discussed in Section \ref{sec:concordances}. The second term consists of Pontrjagin--Weiss classes. As discussed in Section \ref{sec:PWclasses} the lowest of these---in degree $(2n-2)$---corresponds to the configuration space integral associated to the $\Theta$-graph, and so accounts for the class in this degree found for odd $n \leq 399$ by Watanabe \cite{watanabe1}, and show that such classes exist for all $n$. The conclusion of this discussion is depicted in Figure \ref{fig:2}.

In proving Theorem \ref{thm:KrannichRW} Krannich and I were only attempting to calculate within the indicated range, and with the method we used it is not clear how to establish the ``band'' pattern in higher degrees for odd-dimensional discs too. But it does seem feasible that the method used to prove \ref{thm:KupersRW} could be adapted to the odd-dimensional case (though there are significant hurdles) and I think it very likely that the ``band'' pattern  occurs in this case too.

\subsection{Outlook and speculation}

These two theorems have sufficient detail that one is tempted to propose a structural description of $\pi_*(B\Diff_\partial(D^d))\otimes\bQ$. In fact, it seems better to describe $\pi_*(B\Top(d))\otimes\bQ$. Summarising the structural features of the above results, $\pi_*(B\Top(d))\otimes\bQ$ has
\begin{enumerate}[(i)]
\item classes corresponding to Pontrjagin classes, i.e.\ detected by $B\Top(d) \to B\Top$, in degrees $\geq 0$,\label{it:1}
\item classes corresponding to ${K}_{*>0}(\bZ)\otimes\bQ$ if $d$ is odd, in degrees $\gtrsim d$, and a class corresponding to ${K}_0(\bZ)\otimes\bQ$ in degree $d$ if $d$ is even (detected by the Euler class),\label{it:2}
\item classes supported in bands of degrees around $k \cdot d$ for each $k \geq 2$ (at least for $d$ even, but lets suppose that this also occurs for $d$ odd).\label{it:3}
\end{enumerate}

\clearpage

\vspace{3ex}

\noindent\textbf{Orthogonal calculus.}
This behaviour could be explained by Weiss' theory of orthogonal calculus \cite{WeissOrthogonal}, a calculus of functors for continuous functors $\mathsf{F}: \mathcal{J} \to \mathcal{T}op$ defined on the category $\mathcal{J}$ of real inner product spaces and their isometric embeddings. It may be applied to the functor $\mathsf{Bt} : V \mapsto B\Top(V)$, where it provides a tower of Taylor approximations
\vspace{-2ex}
\[
\begin{tikzcd}
 & \vdots\\[-20pt]
&T_2\mathsf{Bt}(V) \dar & \Omega^\infty(S^{2 \cdot V} \wedge \Theta\mathsf{Bt}^{(2)})_{h\OO(2)} \lar\\
&T_1\mathsf{Bt}(V) \dar & \Omega^\infty(S^{1 \cdot V} \wedge \Theta\mathsf{Bt}^{(1)} )_{h\OO(1)} \lar\\
\mathsf{Bt}(V) \arrow[r] \arrow[ru] \arrow[ruu]& T_0\mathsf{Bt}(V) \arrow[r, equal] & B\Top,
\end{tikzcd}
\]
whose $k$th layer is described in terms of an $\OO(k)$-spectrum $\Theta\mathsf{Bt}^{(k)}$, the \emph{$k$-th derivative}, and the 1-point compactifications $S^{k \cdot V}$ of the vector spaces $\bR^k \otimes V$. The zeroth Taylor approximation is the stabilisation of the functor, in this case $B\Top$. It is a theorem of Waldhausen that $\Theta\mathsf{Bt}^{(1)}$ is $\mathrm{A}(*) = \mathrm{K}(\bS)$, with a certain $\OO(1)$-action, and in view of the rational equivalence $\mathrm{K}(\bS) \to \mathrm{K}(\bZ)$ points \ref{it:1} and \ref{it:2} above can be accounted for by the first Taylor approximation. Point \ref{it:3} would then be accounted for if
\begin{enumerate}[(i)]
\item the Taylor tower converges (rationally), and

\item for each $k \geq 2$ the homotopy orbits $(\Theta\mathsf{Bt}^{(k)})_{h\SO(k)}$ of the derivative spectra have finitely-many nontrivial (rational) homotopy groups .
\end{enumerate}
In this worldview, the (finitely-many) rational homotopy groups of $(\Theta\mathsf{Bt}^{(k)})_{h\SO(k)}$ correspond to the rational homotopy classes of $B\Top(d)$ in the $k$-th band which are not detected by Pontrjagin classes or algebraic $K$-theory: more precisely, the residual $\OO(k)/\SO(k)$-action splits $\pi_*((\Theta\mathsf{Bt}^{(k)})_{h\SO(k)})\otimes\bQ$ into eigenspaces, and the $(-1)^d$-eigenspace provides the $k$th band of $B\Top(d)$. In particular, this worldview predicts that the homotopy groups in the $k$th band depend only on the parity of $d$.

\vspace{2ex}

By Theorems \ref{thm:KupersRW} and \ref{thm:KrannichRW} this property does indeed hold for the second band, and Krannich and I \cite{KrannichR-WOdd} have used this to investigate the second derivative $\Theta\mathsf{Bt}^{(2)}$, establishing a rational equivalence
$$\Theta\mathsf{Bt}^{(2)} \simeq_\bQ \mathrm{map}(S^1_+, \bS^{-1})$$
of $\OO(2)$-spectra, where $\OO(2)$ acts in the usual way on $S^1$. Reis and Weiss \cite{ReisWeiss2016} had earlier shown that $\mathrm{map}(S^1_+, \bS^{-1})$ is the second derivative of the orthogonal functor $\mathsf{Bg}(V) := B\mathrm{hAut}(S(V))$, the classifying space of the monoid of homotopy automorphisms of the unit sphere in the inner product space $V$, and the natural map $\mathsf{Bt} \to \mathsf{Bg}$ (in fact zig-zag) induces an equivalence on rationalised second derivatives.

\vspace{3ex}

\noindent\textbf{Automorphisms of little discs operads and graph complexes.}
In Section \ref{sec:CSI} I explained that there is a map
\begin{equation}\label{eq:AutOpMap}
B\Top(d) \lra B\hAut(E_d)
\end{equation}
corresponding to a derived action of $\Top(d)$ on the little $d$-discs operad $E_d$. I mentioned also that the derived automorphisms of the rationalisation $E_d^\bQ$ have been analysed by Fresse, Turchin, and Willwacher \cite{FTW}, giving an identification $\pi_i(\hAut(E_d^\bQ)) = H_i(\mathrm{GC}_{d}^2)$ for $d \geq 3$ in terms of a version of Kontsevich's graph complex. There is a loop-order decomposition $\mathrm{GC}_{d}^2 = \bigoplus_{g \geq 1} \mathrm{GC}_{d}^{2, g\text{-loop}}$ and, for $g \geq 2$,
$$H_*(\mathrm{GC}_{d}^{2, g\text{-loop}}) \text{ is supported in degrees } * \in [g(d-3)+3,g(d-2)+1],$$
and furthermore up to translating degrees this homology depends only on the parity of $d$. There are some computer calculations of these groups, but they are largely unknown. On the other hand the 1-loop part is completely known, and is
$$H_*(\mathrm{GC}_{d}^{2, 1\text{-loop}}) = \bigoplus_{\mathclap{\substack{k \geq 1, \\ k \equiv 2d+1 \mod 4}}} \bQ[d-k].$$

Writing $E_V$ for the little discs operad modelled on the unit disc in an inner product space $V$, one can consider orthogonal calculus applied to the functor $\mathsf{Ba} : V \mapsto B\hAut(E_V)$ (or perhaps better $B\hAut(E_V^\bQ)$: there is an important and subtle question of whether $B\hAut(E_V) \to B\hAut(E_V^\bQ)$ is a rationalisation on universal covers, which I shall elide). Presumably by passing to appropriate models one can upgrade the maps \eqref{eq:AutOpMap} to a map $\mathsf{Bt} \to \mathsf{Ba}$ of orthogonal functors. The data above would seem to suggest that $\mathsf{Ba}$ enjoys precisely the property described in the last section: that $\Theta\mathsf{Ba}^{(1)}$ is rationally equivalent to the $\OO(1)$-spectrum $\map(\mathbb{CP}^\infty_+, \bS^{0})$, where the action is by complex conjugation, and that for $k \geq 2$ the rational homotopy groups of $(\Theta\mathsf{Ba}^{(k)})_{h\SO(k)}$ are supported in degrees $[4-3k,2-2k]$ and combine the $k$-loop graph homology for both parities encoded as the $\OO(k)/\SO(k)$-eigenspace decomposition. (In particular this suggests an action of the twisted group ring $H^{-*}(B\SO(k);\bQ)[\OO(k)/\SO(k)]$ on the graded vector space
$$s^{-2k}H_*(\mathrm{GC}_{2}^{2, k\text{-loop}}) \oplus s^{-3k}H_*(\mathrm{GC}_{3}^{2, k\text{-loop}}),$$
giving a potentially nontrivial relationship between even and odd graph homology.) 

There are two reasons $B\Top(d) \to B\hAut(E_d^\bQ)$ cannot be a rational equivalence:
\begin{enumerate}[(i)]

\item\label{it:Non1} this map tends to kill the Pontrjagin--Weiss classes (and for $d$ odd the algebraic $K$-theory classes),

\item\label{it:Non3} the 1-loop graph contribution $H_{*+1}(\mathrm{GC}_d^{2, 1\text{-loop}}) \subset \pi_*(B\hAut(E_d^\bQ))$ does not come from $\pi_*(B\Top(d))\otimes\bQ$.

\end{enumerate}

Point \ref{it:Non1} concerns the contribution to $B\Top(d) = \mathsf{Bt}(\bR^d)$ due to the first Taylor approximation $T_1 \mathsf{Bt}(\bR^d)$, and, given the degrees in which the 1-loop graphs contribute, point \ref{it:Non3} presumably concerns the contribution to $B\hAut(E_d^\bQ) = \mathsf{Ba}(\bR^d)$ due to the first Taylor approximation $T_1 \mathsf{Ba}(\bR^d)$. Together these suggest that a better question is to ask about the rational homotopy cartesianness of
\begin{equation}\label{eq:RefinedQ}
\begin{tikzcd}
B\Top(d) \arrow[r,equals] & \mathsf{Bt}(\bR^d) \dar \rar& T_1 \mathsf{Bt}(\bR^d) \dar\\
B\hAut(E_d^\bQ) \arrow[r,equals] & \mathsf{Ba}(\bR^d) \rar& T_1 \mathsf{Ba}(\bR^d).
\end{tikzcd}
\end{equation}
If this square were rationally homotopy cartesian for all large enough $d$ then in particular the maps on derivatives $\Theta\mathsf{Bt}^{(k)} \to \Theta\mathsf{Ba}^{(k)}$ would be rational equivalences for all $k \geq 2$. As evidence for this, from the proposed description of the rational homotopy groups of $(\Theta\mathsf{Ba}^{(k)})_{h\SO(k)}$ described above, and the calculation of 2-loop graph homology, one can easily deduce that $\Theta\mathsf{Ba}^{(2)} \simeq_\bQ \mathrm{map}(S^1_+, \bS^{-1})$ and that the induced map $\Theta\mathsf{Bt}^{(2)} \to \Theta\mathsf{Ba}^{(2)}$ is indeed a rational equivalence.

When $d=2n$, rational homotopy cartesianness of \eqref{eq:RefinedQ} is equivalent to the map
$$B\Diff_\partial^\fr(D^{2n}) \simeq \Omega^{2n}_0 \Top(2n) \lra \Omega^{2n}_0(\hAut(E_{2n}^\bQ) \times \Top)$$
being a rational equivalence, and by comparing certain graphical models arising in \cite{KR-WDisks} with graph complexes arising in the operadic model for embedding calculus it looks like this might be true. Kupers, Willwacher, and I are trying to make this precise.

\section{Methods}

I will now explain some of the ideas which go into the proofs of Theorems \ref{thm:KupersRW} and \ref{thm:KrannichRW}, though my goal is to give an overall impression of the methods involved rather than explain how exactly they are combined to prove these two particular results.

\subsection{Weiss fibre sequences and the general strategy}\label{sec:WeissSeq}

Weiss' proof of the existence of Pontrjagin--Weiss classes contains an observation  \cite[Remark 2.1.3]{WeissDalian} which technically does not play a role in his argument but is central to its philosophy. It is more a general principle than a specific formulation, and I shall give it only in a modestly general form: variants of it underlie many recent results about diffeomorphism groups \cite{kupersdisk, KR-WAlg, KrannichConc, KnudsenKupers, KR-WDisks, BustamanteRW, BKK, KrannichR-WOdd}. Let $W$ be a manifold with boundary $\partial W$ decomposed as $\partial_- W \cup \partial_+ W$ into codimension zero submanifolds with common boundary. Then there is a fibration
\begin{equation}\label{eq:WeissFibSeq}
B\Diff_\partial(\partial_- W \times [0,1]) \lra B\Diff_\partial(W) \lra B\Emb_{\partial_+ W}^{\cong}(W, W).
\end{equation}
The rightmost term needs a little explaining: it is the classifying space of the group-like topological monoid $\Emb_{\partial_+ W}^{\cong}(W, W)$ of those self-embeddings of $W$ which are the identity on $\partial_+ W$, and which are isotopic to diffeomorphisms. But crucially these self-embeddings are allowed to send $\partial_- W$ into the interior of $W$, as indicated in the figure below.

\vspace{-2ex}

\begin{figure}[h]
\begin{center}
\includegraphics[width=8cm]{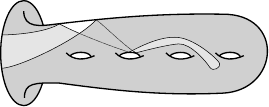}
\end{center}
\end{figure}

\vspace{-2ex}

This is not a technically difficult result (after passing to a different model of the rightmost term it is a simple consequence of the parameterised isotopy extension theorem). Somewhat more technical is Kupers' theorem \cite[Section 4]{kupersdisk} that this fibration sequence deloops---with respect to the evident composition law on $B\Diff_\partial(\partial_- W \times [0,1])$---which is sometimes convenient. 

The importance of this fibration sequence is the following strategy which it indicates: to understand $B\Diff_\partial(\partial_- W \times [0,1])$, you can instead try to understand $B\Diff_\partial(W)$ and $B\Emb_{\partial_+ W}^{\cong}(W, W)$, \emph{for any manifold $W$ containing $\partial_- W$ in its boundary}. This is powerful because these two spaces can sometimes be accessed, though by very different methods: the homology of $B\Diff_\partial(W)$ by parameterised surgery theory, and the homotopy of $B\Emb_{\partial_+ W}^{\cong}(W, W)$ by embedding calculus. Let me explain how this strategy may be implemented to study $B\Diff_\partial(D^d)$.

\vspace{3ex}

\noindent\textbf{Set-up for discs.} If $d=2n$ then take $W_{g,1} := D^{2n} \# g(S^n \times S^n)$ and $\partial_-W_{g,1} = D^{2n-1} \subset \partial W_{g,1}$, so that the fibration \eqref{eq:WeissFibSeq} takes the form
\begin{equation}\label{eq:WeissSeqEven}
B\Diff_\partial(D^{2n}) \lra B\Diff_\partial(W_{g,1}) \lra B\Emb_{\partial_+ W_{g,1}}^{\cong}(W_{g,1}, W_{g,1}),
\end{equation}
where we have identified $D^{2n-1} \times [0,1] \approx D^{2n}$.

If $d=2n+1$ then instead of directly using \eqref{eq:WeissFibSeq}, take the handlebody $V_g := \natural g (S^n \times D^{n+1})$, and use a variant of \eqref{eq:WeissFibSeq} allowing most of the boundary of $V_g$ to not be fixed. This takes the form
\begin{equation}\label{eq:WeissSeqOdd}
B\mathrm{C}(D^{2n}) \lra B\Diff_{D^{2n}}(V_g) \lra B\Emb^{\cong}_{D^{2n}, W_{g,1}}(V_g, V_g),
\end{equation}
where the left-hand term is given by pseudoisotopies of $D^{2n}$, the middle term by diffeomorphisms of $V_g$ fixing a disc $D^{2n} \subset \partial V_g$ but allowing the rest of the boundary to move, and the right-hand term is given by self-embeddings of $V_g$ which preserve $\partial V_g \setminus \mathrm{int}(D^{2n}) = W_{g,1}$ setwise, and furthermore preserve a disc $D^{2n} \subset W_{g,1}$ pointwise. By the fibration sequence 
\begin{equation}\label{eq:DiscsConc}
B\Diff_\partial(D^{2n+1}) \lra B\mathrm{C}(D^{2n}) \lra B\Diff_\partial(D^{2n})\end{equation}
from \eqref{eq:ConcFibSeq}, given $B\Diff_\partial(D^{2n})$ it is equivalent to get at $B\Diff_\partial(D^{2n+1})$ or $B\mathrm{C}(D^{2n})$, and I will explain below why the latter is more accessible.

\vspace{3ex}

\noindent\textbf{Parameterised surgery.} 
 The reason for the choice of manifold $W_{g,1}$ is that Galatius and I \cite{grwstab1} have shown that the maps $B\Diff_\partial(W_{g,1}) \to B\Diff_\partial(W_{g+1,1})$, induced by the evident embeddings $W_{g,1} \hookrightarrow W_{g+1,1}$, are homology isomorphisms in a range of homological degrees tending to infinity with $g$, as long as $2n \geq 6$, and furthermore \cite{grwcob} that a certain parameterised Pontrjagin--Thom map
$$\hocolim_{g \to \infty} B\Diff_\partial(W_{g,1}) \lra \Omega^\infty_0 \mathrm{MT}\theta^n,$$
to the infinite loop space of a certain Thom spectrum, induces an isomorphism on homology. The rational cohomology of the right-hand side is quite simple: it is a polynomial algebra on certain easily-defined cohomology classes, known as Miller--Morita--Mumford classes.

These results are analogues in high dimensions of Harer's \cite{Ha} theorem on the stability of the homology of mapping class groups of oriented surfaces, and Madsen and Weiss' \cite{MW} theorem on the stable homology of these mapping class groups. In fact, the stability result holds much more generally for all $2n$-manifolds of the form $W \# g(S^n \times S^n)$ with $2n \geq 6$ and $W$ simply-connected (or even with virtually polycyclic fundamental group \cite{Friedrich}), and there is an analogous description of the stable homology for any $W$ of any even dimension \cite{grwstab2} (including dimension 4). See Galatius' contribution to the 2014 ICM for an overview of this theory. 

In odd dimensions the stable homology of the diffeomorphism groups of the analogous manifolds $D^{2n+1} \# g(S^n \times S^{n+1})$ is not yet known, but Botvinnik and Perlmutter \cite{perlmutterbotvinnik} have a version for $B\Diff_{D^{2n}}(V_g)$, and Perlmutter \cite{perlmutterstab} has the appropriate stability theorem in this case. This accounts for the use of the modified Weiss fibre sequence \eqref{eq:WeissSeqOdd} in odd dimensions, rather than a more obvious analogue of \eqref{eq:WeissSeqEven} involving $D^{2n+1} \# g(S^n \times S^{n+1})$.

\vspace{3ex}

\noindent\textbf{Embedding calculus.} The difficulty of studying embeddings of one manifold into another depends on the codimension, but this must be counted appropriately. What matters is the \emph{geometric dimension} of the target minus the \emph{handle dimension} of the source. In particular, if $W$ is $d$-dimensional but can be constructed from $\partial_+ W \times [0,1]$ by attaching handles of index $\leq h$, then self-embeddings of $W$ relative to $\partial_+ W$ have codimension $d-h$. If this codimension is $\geq 3$, then the theory of embedding calculus as developed by Goodwillie, Klein, and Weiss \cite{weissembeddings, goodwillieweiss, goodwillieklein} can be used to access spaces of self-embeddings of $W$ relative to $\partial_+ W$. This theory provides a tower
\vspace{-1ex}
\begin{equation}\label{eq:EmbTower}
\begin{tikzcd}[column sep=0.5cm]
 & \vdots\\[-20pt]
& BT_3\Emb_{\partial_+ W}^{\cong}(W,W)\dar \\
& BT_2\Emb_{\partial_+ W}^{\cong}(W,W) \dar \\
B\Emb_{\partial_+ W}^{\cong} (W,W) \arrow[r] \arrow[ru] \arrow[ruu]& BT_1\Emb_{\partial_+ W}^{\cong}(W,W),
\end{tikzcd}
\end{equation}
such that as long as the codimension (as described above) is $\geq 3$ the map
$$B\Emb_{\partial_+ W}^{\cong} (W,W) \lra B T_\infty\Emb_{\partial_+ W}^{\cong} (W,W) = \holim_{k \to \infty} B T_k\Emb_{\partial_+ W}^{\cong} (W,W)$$
is an equivalence. The bottom stage $BT_1\Emb_{\partial_+ W}^{\cong}(W,W)$ is equivalent to the classifying space of the monoid $\Bun_{\partial_+ W}^{\cong}(TW, TW)$ of bundle maps $TW \to TW$ which are the identity over $\partial_+ W$ and which are homotopic to the derivative of a diffeomorphism, and the homotopy fibre of $BT_k \Emb_{\partial_+ W}^{\cong}(W,W) \to BT_{k-1} \Emb_{\partial_+ W}^{\cong}(W,W)$ has a description in terms of a space of sections of a bundle $Z_k \to C_k(W)$ over the configuration space of $k$ unordered points in $W$, whose fibres are constructed from the configuration spaces of $\leq k$ ordered points in $W$. Thus in principle the bottom stage and these homotopy fibres are amenable to calculation by homotopical methods.

For the manifold $W_{g,1}$ and $\partial_+ W_{g,1} = D^{2n}$, the codimension in the sense described is $n$, so the embedding calculus tower converges as long as $2n \geq 6$. For the manifolds $V_g$ there is a similar tower for $B\Emb^{\cong}_{D^{2n}, W_{g,1}}(V_g, V_g)$, which converges for $2n+1 \geq 7$.

\vspace{1ex}

The strategy which suggests itself is then to calculate as much as you can about the middle and right-hand terms of \eqref{eq:WeissSeqEven} and \eqref{eq:WeissSeqOdd} using these two very different methods, and then use these fibre sequences and \eqref{eq:DiscsConc} to deduce things about $B\Diff_\partial(D^d)$. This is a very attractive picture, but for getting explicit answers there is a serious 

\vspace{3ex}

\noindent\textbf{Difficulty.} Parameterised surgery fundamentally gets at the \emph{homology} groups of diffeomorphism groups, whereas embedding calculus, at least if applied in the most classical way, naturally allows one to get at the \emph{homotopy} groups of embedding spaces.

\subsection{Qualitative results}

One situation in which this Difficulty is not so serious is if one wishes to obtain qualitative results about $B\Diff_\partial(\partial_- W \times [0,1])$, for example that its homology or homotopy groups lie in a given Serre class. This was pioneered by Kupers \cite{kupersdisk}, to prove that the homotopy (or equivalently, homology) groups of $B\Diff_\partial(D^d)$ are finitely-generated for $d \neq 4,5,7$. A slight variant of his line of reasoning is as follows. 

Firstly, for $D^{2n}$ consider the Weiss fibre sequence \eqref{eq:WeissSeqEven}, which may be delooped. Using \cite{grwcob, grwstab1}, as long as $2n \geq 6$ the homology of $B\Diff_\partial(W_{g,1})$ is easily seen to be finitely-generated in degrees $* \leq \tfrac{g-3}{2}$, so it suffices to show that the homology of $X := B\Emb_{\partial_+ W_{g,1}}^{\cong}(W_{g,1}, W_{g,1})$ is finitely-generated too. Using the embedding calculus tower \eqref{eq:EmbTower} it is not difficult to show that the higher homotopy groups of $X$ are all finitely-generated, and hence to deduce that the homology of the universal cover $\widetilde{X}$ is finitely-generated. It remains to study the spectral sequence for the fibration
$$\widetilde{X} \lra X \lra B\pi_1(X)$$
and the crucial point here is that the group 
$$\pi_1(X) = \pi_0(\Emb_{\partial_+ W_{g,1}}^{\cong}(W_{g,1}, W_{g,1})) \cong \pi_0(\Diff_\partial(W_{g,1}))/\pi_0(\Diff_\partial(D^{2n}))$$
enjoys Wall's finiteness property ($F_\infty$). In this case it is clear by Kreck's \cite{kreckisotopy} calculation of the group $\pi_0(\Diff_\partial(W_{g,1}))$, but as a general principle it follows from Sullivan's theorem \cite{Su3} that mapping class groups of simply-connected manifolds of dimension $\geq 5$ are commensurable (up to finite kernel, see \cite{KrannichRW}) to arithmetic groups.

Secondly, for $D^{2n+1}$ it suffices, given the above, to prove finite-generation of the homology of $B\mathrm{C}(D^{2n})$, so consider the Weiss fibre sequence \eqref{eq:WeissSeqOdd}, which may also be delooped. Using \cite{perlmutterstab, perlmutterbotvinnik}, as long as $2n+1 \geq 9$ the homology of $B\Diff_{D^{2n}}(V_g)$ is finitely-generated in a stable range, and embedding calculus considerations as above show that the homology of $B\Emb^{\cong}_{D^{2n}, W_{g,1}}(V_g, V_g)$ is finitely-generated too.

\vspace{1ex}

When working modulo a Serre class one can sometimes also determine the lowest nonvanishing term. Bustamante and I \cite{BustamanteRW} have used the above strategy with a Weiss fibre sequence for the manifolds $X_g := S^1 \times D^{2n-1} \# g(S^n \times S^n)$ to analyse $\pi_*(B\Diff_\partial(S^1 \times D^{2n-1}))_{(p)}$ for $2n \geq 6$ modulo the Serre class of finitely-generated $\bZ_{(p)}$-modules, where we show that it vanishes in degrees $* < \min(2p-3, n-2)$ and is $\bigoplus^\infty \bZ/p$ in degree $2p-3$ as long as $2p-3 < n-2$. This was known in the pseudoisotopy stable range using algebraic $K$-theory methods \cite{GKM}, but our work gives a rather different perspective on this infinitely-generated subgroup.

\vspace{3ex}

\noindent\textbf{Turning the Weiss fibre sequence around.} A further point of view on the Weiss fibre sequence is that, assuming $B\Emb_{\partial_+ W}^{\cong}(W, W)$ may be understood using embedding calculus, it reduces questions about $B\Diff_\partial(W)$ for a whole class of manifolds $W$ to questions about the single space $B\Diff_\partial(\partial_- W \times [0,1])$. As the minimal choice of $\partial_-W$ is $\partial_-W = D^{d-1}$, this gives another reason to be particularly interested in $B\Diff_\partial(D^d)$.

Kupers \cite{kupersdisk} exploits this point of view to show---given the homological, and so also homotopical, finite generation of $B\Diff_\partial(D^d)$ discussed above---that $B\Diff_\partial(W)$ has finitely-generated higher homotopy groups for \emph{any} closed 2-connected manifold $W$ of dimension $d \neq 4,5,7$. More recently Bustamante, Krannich, and Kupers \cite{BKK} have extended this to any closed manifold of dimension $2n \geq 6$ with finite fundamental group. 

\subsection{Quantitative results}

To obtain quantitative results one must confront the Difficulty. The most obvious way to do this---in the case of discs---is to try to make Kupers' method from the last section quantitative, by trying to calculate the  homology of the right-hand terms of \eqref{eq:WeissSeqEven} and \eqref{eq:WeissSeqOdd}. 

This is the strategy I pursued with Krannich in \cite{KrannichRW} to prove Theorem \ref{thm:KrannichRW}, and (though for $S^1 \times D^{2n-1}$ and not for discs) with Bustamante in \cite{BustamanteRW}. An important preliminary simplification is to consider the Weiss fibre sequence with framings (or similar): for example, in the framed version 
\begin{equation}\label{eq:WeissSeqEvenFr}
B\Diff_\partial^\fr(D^{2n}) \lra B\Diff_\partial^\fr(W_{g,1}) \lra B\Emb_{\partial_+ W_{g,1}}^{\fr, \cong}(W_{g,1}, W_{g,1})
\end{equation}
of the sequence \eqref{eq:WeissSeqEven}, $B\Diff_\partial^\fr(D^{2n})$ differs from $B\Diff_\partial(D^{2n})$ by a copy of $\Omega^{2n} \OO(2n)$, whose rational homotopy groups are completely understood, but \cite{grwcob, grwstab1} shows that the rational homology of $B\Diff_\partial^\fr(W_{g,1})$ is \emph{trivial} in the stable range, which is far simpler than the rational homology of $B\Diff_\partial(W_{g,1})$. Luckily the effect on the homology of the self-embeddings term is also beneficial. The way we calculate the homology of the (framed) self-embedding spaces in these papers is not in fact using embedding calculus as I have been advertising, but rather using disjunction theory (which is in any case the fuel which makes embedding calculus work \cite{goodwillieweiss, goodwillieklein}, but using it directly is sometimes more convenient), in the form of Morlet's lemma of disjunction in \cite{BustamanteRW} and Goodwillie's multi-relative disjunction lemma \cite{goodwilliethesis} in \cite{KrannichRW}. The nature of the calculations involved makes it hard to say anything very general about them, so I shall not try to.

Instead I should like to discuss an alternative strategy, which is what Kupers and I do in \cite{KR-WDisks} and prepare for in the companion papers \cite{KR-WAlg, KR-WTorelli, KR-WTorelliLie}, and is what leads to the proof of Theorem \ref{thm:KupersRW}. There we adopt the view that embedding calculus is very well suited to calculating---or estimating---the rational \emph{homotopy} groups of $B\Emb_{\partial_+ W_{g,1}}^{\cong}(W_{g,1}, W_{g,1})$, and so we propose to calculate---or estimate---the rational \emph{homotopy} groups of $B\Diff_\partial(W_{g,1})$ (in fact we consider the framed version $B\Diff_\partial^\fr(W_{g,1})$, but again the difference on rational homotopy groups is very mild). Describing these is an interesting problem in its own right, especially in view of Berglund and Madsen's \cite{berglundmadsen2} calculation of the rational homotopy and stable cohomology of the groups of block diffeomorphisms and of homotopy automorphisms of $W_{g,1}$. In the remainder I will focus on this calculation, and not try to explain exactly how it implies Theorem \ref{thm:KupersRW}.

\subsection{Torelli groups}

Diffeomorphisms of $W_{g,1}$ induce automorphisms of $H_n(W_{g,1};\bZ)$ which preserve the intersection form, giving a homomorphism
$$\alpha_g : \Diff_\partial(W_{g,1}) \lra G_g := \begin{cases}
\OO_{g,g}(\bZ) & n \text{ even}\\
\Sp_{2g}(\bZ) & n \text{ odd}.
\end{cases}$$
This is surjective if $n$ is even or $n=1,3,7$, but for other odd values of $n$ has image a certain finite-index subgroup $G'_g$. By analogy with the case $2n=2$, the kernel of $\alpha_g$ is called the \emph{Torelli group} and denoted $\Tor_\partial(W_{g,1})$, so that there is a fibration sequence
\begin{equation}\label{eq:TorelliFib}
B\Tor_\partial(W_{g,1}) \lra B\Diff_\partial(W_{g,1}) \lra BG'_g.
\end{equation}

Now the fundamental group of $B\Diff_\partial(W_{g,1})$ is quite complicated, as it surjects onto the arithmetic group $G'_g$, so although the results of \cite{grwcob, grwstab1} describe the rational cohomology of this space, there is no reason to think that this has much to do with its rational higher homotopy groups. However, by \eqref{eq:TorelliFib} these higher homotopy groups are the same as those of $B\Tor_\partial(W_{g,1})$, and as long as $2n \geq 6$ the Weiss fibre sequence can be used (in ``qualitative mode'') to prove that the space $B\Tor_\partial(W_{g,1})$ is nilpotent \cite[Theorem C]{KR-WAlg}. Thus $B\Tor_\partial(W_{g,1})$ has a meaningful rationalisation, and its rational homotopy and cohomology groups are closely related (in the sense that there are spectral sequences computing each from the other). On the other hand, passing to the infinite covering space $B\Tor_\partial(W_{g,1})$ of $B\Diff_\partial(W_{g,1})$ has an unknown effect on cohomology, so the problem is now to understand the rational cohomology of $B\Tor_\partial(W_{g,1})$.

\vspace{3ex}

\noindent\textbf{Cohomology of Torelli groups.}
The fibration \eqref{eq:TorelliFib} provides a representation of the arithmetic group $G'_g$ on the rational vector spaces $H^i(B\Tor_\partial(W_{g,1});\bQ)$, and as long as $2n \geq 6$ and $g \geq 2$ a further application of the Weiss fibre sequence in ``qualitative mode'' shows \cite[Theorem A]{KR-WAlg} that these are \emph{algebraic representations} of $G'_g$, i.e.\ they extend to representations of the ambient algebraic group $\mathbf{G}_g \in \{\mathbf{O}_{g,g},\mathbf{Sp}_{2g}\}$. As $\mathbf{G}_g$-representations are semisimple, and the irreducibles are classified in terms of Young diagrams and are all summands of tensor powers of the defining $\mathbf{G}_g$-representation $H := H_n(W_{g,1};\bQ)$, the $G'_g$-representation $H^i(B\Tor_\partial(W_{g,1});\bQ)$ may be determined in terms of the vector spaces
\begin{equation}\label{eq:Invariants}
[H^i(B\Tor_\partial(W_{g,1});\bQ) \otimes H^{\otimes S}]^{G'_g}
\end{equation}
for all finite sets $S$, and the structure maps between them given by applying permutations and contractions $H \otimes H \to \bQ$. On the other hand, the vector spaces \eqref{eq:Invariants} are related, by the Serre spectral sequence for \eqref{eq:TorelliFib}, to the cohomology groups
\begin{equation}\label{eq:TwistCoh}
H^*(B\Diff_\partial(W_{g,1}); \mathcal{H}^{\otimes S})
\end{equation}
with coefficients in the $S$-th tensor power of the local system $\mathcal{H}$ on $B\Diff_\partial(W_{g,1})$ provided by the $G'_g$-representation $H$. Using work of Borel \cite{borelstable2} this Serre spectral sequence can be shown to degenerate: Ebert and I \cite{oscarjohannestorelli} introduced this strategy, but as the results of \cite{grwcob, grwstab1} only describe the cohomology of $B\Diff_\partial(W_{g,1})$ with constant coefficients, we were only able to use it to determine $[H^i(B\Tor_\partial(W_{g,1});\bQ)]^{G'_g}$ in a stable range. However, the results of \cite{grwcob, grwstab1} also apply to $B\Diff_\partial^\theta(W_{g,1})$ for quite arbitrary tangential structures $\theta$ (such as framings, but also including ``maps to a space $Y$''),  and exploiting the functoriality of the result with respect to $\theta$ (this kind of argument originates in \cite{R-Wtwist}) it is possible to calculate \eqref{eq:TwistCoh} in a stable range, and hence by the strategy outlined here to calculate $H^*(B\Tor_\partial(W_{g,1});\bQ)$ in a stable range of degrees, as a $\bQ$-algebra and as a $G'_g$-representation. This is done in \cite{KR-WTorelli}.

A similar strategy can be applied to $B\Diff^\fr_\partial(W_{g,1})$, though by a subtlety in the proof the argument above does not directly calculate the cohomology of $B\Tor^\fr_\partial(W_{g,1})$. Instead, there is a certain fibration sequence
$$X_1(g) \lra B\Tor^\fr_\partial(W_{g,1}) \lra X_0$$
with $X_0$ a loop space having rational cohomology $\Lambda_\bQ[\bar{\sigma}_{4j-2n-1} \, | \, j > n/2]$. The $\bar{\sigma}_{4j-2n-1}$ are secondary characteristic classes associated to the fact that the family signature vanishes for two different reasons on $B\Tor^\fr_\partial(W_{g,1})$: because of the framing, and because of the triviality of the action on $H^n(W_{g,1};\bZ)$. The analogue of the argument above leads to the following description of the cohomology of the nilpotent space $X_1(g)$. For $r \geq 3$ and $v_1, \ldots, v_r \in H^n(W_{g,1};\bQ)$ there are defined \emph{twisted Miller--Morita--Mumford classes}
$$\kappa_1(v_1 \otimes \cdots \otimes v_r) \in H^{(r-2)n}(X_1(g);\bQ)$$
which satisfy:
\begin{enumerate}[(i)]
\item linearity in each $v_i$,\label{it:X11}

\item $\kappa_1(v_{\sigma(1)} \otimes \cdots \otimes v_{\sigma(r)}) = \mathrm{sign}(\sigma)^n \cdot \kappa_1(v_1 \otimes \cdots \otimes v_r)$,\label{it:X12}

\item $\sum_i \kappa_1(v_1 \otimes \cdots \otimes v_{k-1} \otimes a_i) \smile \kappa_1(a_i^\# \otimes v_k \otimes \cdots \otimes v_r) = \kappa_1(v_1 \otimes \cdots \otimes v_r)$,\label{it:X13}

\item $\sum_i \kappa_1(v_1 \otimes \cdots \otimes v_r \otimes a_i \otimes a_i^\#)=0$,\label{it:X14}
\end{enumerate}
where $\sum_i a_i \otimes a_i^\# \in H^n(W_{g,1};\bQ)^{\otimes 2}$ is dual to the intersection form. The framed analogue of the group $G'_g$ acts on $\kappa_1(v_1 \otimes \cdots \otimes v_r)$ by its evident action on the $v_i \in H^n(W_{g,1};\bQ)$. Kupers and I \cite{KR-WDisks} show that, in a range of degrees tending to infinity with $g$, the cohomology algebra of $X_1(g)$ is generated by the classes $\kappa_1(v_1 \otimes \cdots \otimes v_r)$ and subject only to the relations \ref{it:X11}--\ref{it:X14}.

\vspace{3ex}

\noindent\textbf{Homotopy of Torelli groups.}
As the rational cohomology of $X_1(g)$ is supported in degrees divisible by $n$ in a stable range, it follows formally that its rational homotopy groups are supported in degrees $\cup_{r \geq 1} [r(n-1)+1,rn]$ in this stable range, so exhibit a band pattern. 

But it turns out that we can do a lot better. It is not hard to see that the above data in fact presents a quadratic algebra, generated by the elements $\kappa_1(v_1 \otimes v_2 \otimes v_3)$ of degree $n$ (modulo \ref{it:X14}), and it is then tempting to ask whether this quadratic algebra is Koszul.  Kupers and I \cite{KR-WTorelliLie} prove that $H^*(X_1(g);\bQ)$ is indeed Koszul in a stable range of degrees (this was simultaneously proved by Felder, Naef, and Willwacher \cite{FNW}), so it follows that in this range $\pi_*(X_1(g))\otimes\bQ$ is in fact supported in degrees of the form $r(n-1)+1$, and is furthermore given by the quadratic dual Lie algebra. Up to two extension questions, this calculates $\pi_*(B\Diff_\partial(W_{g,1}))\otimes\bQ$ in a stable range.

\vspace{2ex}

\noindent\textbf{Acknowledgements.} It is a pleasure to thank those with whom I have collaborated on the ideas described here: Mauricio Bustamante, Johannes Ebert, S{\o}ren Galatius, Manuel Krannich, and Alexander Kupers.

\bibliographystyle{amsplain}
\bibliography{biblio}

\end{document}